\documentclass[12pt,amsfonts,amssymb]{article}
\usepackage{amsfonts}
\usepackage{amssymb}
\def\restriction{\upharpoonright}
\newtheorem{theorem}{Theorem}[section]

\newtheorem{lemma}[theorem]{Lemma}
\newtheorem{definition}[theorem]{Definition}
\newtheorem{claim}[theorem]{Claim}
\newtheorem{fact}[theorem]{Fact}
\newtheorem{remark}[theorem]{Remark}

\newtheorem{notation}[theorem]{Notation}
\newtheorem{conjecture}[theorem]{Conjecture}

\newcommand{\qed}{\hfill\mbox{$\Box$}}

\def\eq{\rm eq}

\def\b1K{\mbox{\boldmath $K$}_{-1}}
\def\bK{\mbox{\boldmath $K$}}

\def\Mscr{{\cal M }}

\def\Pscr{{\cal P }}

\def\cf{\mathop{\rm cf}}

\def\th{\mathop{\rm Th}}

\def\stp{\mathop{\rm stp}}

\def\aut{\rm aut}

\def\acl{{\rm acl}}

\def\grpf #1 #2{{\rm grp}_{#2}(#1)}
\def\spanf #1 #2{{\rm span}_{#2}(#1)}
\def\fldf #1 #2{{\rm fld}_{#2}(#1)}
\def\dclf #1 #2{{\rm dcl}_{#2}(#1)}
\def\rclf #1 #2{{\rm rcl}_{#2}(#1)}
\def\aclf #1 #2{{\rm acl}_{#2}(#1)}
\def\acff #1 #2{{\rm acf}_{#2}(#1)}
\def\strf #1 #2{{\rm str}_{#2}(#1)}
\def\tclf #1 #2{{\rm acf}_{#2}(#1)}
\def\abar{\mbox{\boldmath $a$}}
\def\bbar{{\bf b}}
\def\cbar{{\bf c}}
\def\dbar{{\bf d}}

\def\fbar{{\bf f}}

\def\hbar{{\bf h}}

\def\xbar{{\bf x}}

 \def\restrict{|}
\def\tp{{\rm tp }}

\def\id{\rm id}
\def\eq{\rm eq}
\newcommand{\sidebar}[1]{\vskip10pt\noindent
 \hskip.70truein\vrule width2.0pt\hskip.5em
 \vbox{\hsize= 4truein\noindent\footnotesize\relax #1 }\vskip10pt\noindent}
\title{Subsets of superstable structures are weakly benign}
\author{Bektur Baizhanov\\Institute for Informatics and Control\\
Almaty, Kazakhstan\thanks{ Partially supported by  CDRF grant
KM2-2246.}\\
John T. Baldwin\\
Department of Mathematics, Statistics and  Computer
Science\\
University of Illinois at Chicago\thanks{ Partially supported by
NSF grant  DMS-0100594 and CDRF grant KM2-2246.} \and Saharon
Shelah\thanks{This is    Publication 815 in Shelah's bibliography.
This research was  partially supported by The Israel Science
Foundation.
}\\ Department of Mathematics\\
Hebrew University of Jerusalem\\Rutgers University}

\begin{document}
\maketitle
%\sidebar{These are notes on the Shelah approach to
%various conjectures of Baizhanov and Baldwin. I take liberties to
%rewrite in my own notation and language.  These are based on notes
%written by Shelah, Sept. 24/25 2001 and given the number F529. I
%have received considerable further explanation from Shelah. Pages
%0.0 through 0.7 are concerned with ways to rephrase the concept
%$A$ is closed in $C^{\eq}$ and we add `a predicate' to name $A$.
%(Which is manifestly strange since $A$ is in many sorts.) This
%material hasn't been typed up yet and may not be needed in view of
%the introduction of the concept: weakly benign. The next few
%definitions  are in my view nicer ways to describe the concepts
%described on pages 1.1, 1.2, of the notes. }

Baizhanov and Baldwin \cite{BaBa} introduce the notions of benign
and weakly benign sets to investigate the preservation of
stability by naming arbitrary subsets of a stable structure. They
connect the notion with work of Baldwin, Benedikt, Bouscaren,
Casanovas, Poizat, and Ziegler.  Stimulated by \cite{BaBa}, we
investigate here the existence of benign or weakly benign sets.

\begin{definition}
\begin{enumerate}
\item The set $A$ is {\em benign} in $M$ if for every
$\alpha, \beta \in M$ if $p=\tp(\alpha/A) = \tp(\beta/A)$ then
$\tp_*(\alpha/A) = \tp_*(\beta/A)$ where the $*$-type is the type
in the language $L^*$ with a new predicate $P$ denoting $A$.
\item The set $A$ is {\em weakly benign} in $M$ if for every
$\alpha, \beta \in M$ if $p=\stp(\alpha/A) = \stp(\beta/A)$ then
$\tp_*(\alpha/A) = \tp_*(\beta/A)$ where the $*$-type is the type in
language with a new predicate $P$ denoting $A$.
\end{enumerate}
\end{definition}
%\sidebar{Shelah originally wrote $(M,A)$ is stable for $A$ is {\em benign} in $M$
%to indicate a relationship to stability over a predicate. }
\begin{conjecture}[too optimistic]
\label{mostgeneral}  If $M$ is a model of stable theory $T$ and $A \subseteq M$
then $A$ is benign.
\end{conjecture}
Shelah observed, after learning of the Baizhanov--Baldwin
reductions of the problem to equivalence relations, the following
counterexample.
\begin{lemma}  There is an $\omega$-stable rank 2 theory $T$ with ndop
which has a model $M$ and set $A$ such that
$A$ is not { benign} in $M$.
\end{lemma}
Proof:  The universe of $M$ is partitioned into two sets denoted
by $Q$ and $R$.  Let $Q$ denote $\omega \times \omega$   and $R$
denote $\{0,1\}$.  Define $E(x,y,0)$ to hold if the first
coordinates of $x$ and $y$ are the same and $E(x,y,1)$ to hold if
the second coordinates of $x$ and $y$ are the same.  Let $A$
consist of one element from each $E(x,y,0)$-class and one element
of all but one $E(x,y,1)$-class such that no two members of $A$
are equivalent for either equivalence relation.  It is easy to
check that letting $\alpha$ and $\beta$ denote the two elements of
$R$, we have a counterexample. In this case, the type $p$ is
algebraic.  Algebraicity is a completely artificial restriction.
Replace each $\alpha$ and $\beta$ by an infinite set of points
which behave exactly as $\alpha$, $\beta$ respectively.  We still
have a counterexample.  In either case, $\alpha$ and $\beta$ have
different strong types.  This leads to the following weakening of
the conjecture.

%\sidebar{This sidebar begins a different phrasing of
%the weakened conjecture than the one currently favored.  If all
%goes well the notion of `closed' will just disappear.
%\begin{definition}
%We say $A$ is closed if $\acl^{\eq}(A) = A$.
%\end{definition}
%Then a weakening of the conjecture is to assume $A$ is closed.
%Closed is problematic because of the awkwardness of expanding by
%predicates in many sorts.
%One possibility for proving a theorem in the home sort is to
%define weakly closed as: all types over $A$ are stationary, and try
%to prove the conjecture for weakly closed sets.  Shelah's note
%of Dec. 6, 01 has a further idea for working in $C^{\eq}$.   For now,
%we omit that possibility.}
\begin{conjecture}[Revised]
\label{closedgeneral}%\begin{enumerate}
%\item If $M$ is a model of stable theory $T$ and
%$A \subseteq M^{\eq}$
%is closed
%then $A$ is benign.
%\item
If $M$ is a model of stable theory $T$ and
$A$ is an arbitrary subset of $M$
then $A$ is weakly benign.
%\end{enumerate}
\end{conjecture}
We give here a proof of  Conjecture~\ref{closedgeneral} in the
superstable case.  There are two steps. In the first we show that
if $(M,A) $ is not (weakly) benign then there is a certain
configuration within $M$. (This uses only $T$ stable.)  The second
shows that this configuration is contradicted for superstable $T$.
Note that if $(M,A)$ is not weakly benign, neither is any
$L^*$-elementary extension of $(M,A)$ so we may assume any
counterexample is sufficiently saturated.

\section{Refining a counterexample}
In this section we choose a specific way in which sufficiently
saturated pair $(M,A)$ where $\th(M)$ is  stable, fails to be
weakly benign. Fix $(M,A)$, a $\kappa^+$-saturated of a stable
theory $T$ where $\kappa = \kappa^{|T|}$ is regular.

 We introduce some notation.  Recall that $A$ is {\em relatively
 $\kappa$-saturated} in $M$ if every type over (a subset of $A$)
 whose domain has cardinality less than $\kappa$ and which is
 realized in $M$, is also realized in $A$.
First note that for any $c\in M-A$, there is a pair $(M_1,A_1)$
such that $A_1$ is relatively $\kappa$-saturated in $ A$; $A_1
\cup {c} \subseteq M_1$ and $M_1$ is independent from $A$ over
$A_1$; $A_1$ and $M_1$ have cardinality $\kappa$
 and $M_1$ is $\kappa$-saturated.  For this, choose $A_0\subset A$ with $c$
independent from $A$ over $A_0$ and $|A_0|< \kappa$ (which follows
since $\kappa \geq |T| \geq \kappa(T)$). Then extend $A_0$ to a
subset $A_1$ of $A$ with cardinality at most $\kappa$ which is
relatively $\kappa$-saturated in $A$.   Finally,  let $M_1\prec M$
be $\kappa$-prime over $A_0 \cup c$.  We have shown the following
class $\bK_c$ is not empty.
\begin{notation}
\label{nota}
\begin{enumerate}
\item For any $c \in M$, let $\bK_c$ be the class of pairs
 $(M_1,A_1)$ with $c \in M_1 \prec M$ such that
$A_1$ is relatively $\kappa$-saturated in $ A$; $A_1 \cup {c} \subseteq M_1$ and $M_1$ is independent from $A$ over
$A_1$; $A_1$ and $M_1$ have cardinality $\kappa$
 and $M_1$ is $\kappa$-saturated with $|M_1| \leq \kappa$.
\item
For any $a,b$ in $M$ which realize the same type over $A$, let
$\bK^1_{a,b}$ be the set of tuples $\langle A_1, M_a, M_b, N_a,
g\rangle$ such that $(M_a,A_1)$ and $(M_b, A_1)$ are in $\bK_a,
\bK_b$ respectively, $g$ is an isomorphism between $M_a$ and $M_b$
(subsets of $M$) over $A_1$ (taking $a$ to $b$), $N_a$ contains
$M_a$ and  is saturated with cardinality $\kappa$, and $N_a$ is
independent from $A$ over $A_1$.
%$(N_a, N_a \cap A) \prec (M,A)$.
\item  Let $\bK^2_{a,b}$ be the set of tuples $\langle A_1, M_a, M_b, N_a,
g\rangle \in \bK^1_{a,b}$ such that $g$ is  an isomorphism
 between $M^{\eq}_a$ and $M^{\eq}_b$ over $A^{\eq}_1$.
 \item We will write $K^i$ to denote either $K^1$ or $K^2$.  Note
 the only difference between them is that $K^2$ has a more
 restrictive requirement on the isomorphism $g$.
\end{enumerate}
\end{notation}

Note that the last clause of item 2 implies that $N_a$ is
independent from $A$ over $N_a \cap A$ and that $ N_a  \cap A =
A_1 = M_a \cap A $.
%\sidebar{I don't
%understand the very last clause of item 2 or where it is used
%either here or in section 2. Note that the last clause of item 2
%implies that $N_a$ is independent from $A$ over $N_a \cap A$. But
%I don't see where that is used. }
%\sidebar{ this strengthen (2) in 0.9, so will help in the proof for superstable
%hopefully elsewhere   ss}
%\sidebar{remark helps but still need to think a bit. jtb}
Moreover, if $\langle A_1, M_a, M_b, N_a, g\rangle \in \bK_{a,b}$
and $B \subseteq A $ with $|B|\leq \kappa$ then there is an
$\langle A'_1, M'_a, M'_b,N'_a, g' \rangle\in \bK_{a,b}$ with $A_1
\cup B \subseteq A'_1$.  (Just include $B$ when making the
construction from the first paragraph of this section to show
$\bK_{a,b}$ is nonempty).  We need a couple of other properties of
$\bK_{a,b}$. Note that $\bK_{a,b}$ is naturally partially ordered
by coordinate by coordinate inclusion.

\begin{lemma}\label{chains} Every
increasing chain from $\bK^i_{a,b}$ of length $\delta$ a limit
ordinal less than $\kappa^+$ has an upper bound in $\bK^i_{a,b}$.
 \end{lemma}

 Proof. If the cofinality of
the chain is at least $\kappa_r(T)$,  just take the union (in each
coordinate).  We check that  $ N^\delta_a$, $  A$ are independent
 over  $ A^\delta$:  By induction, for every
$\alpha<\beta<\delta$, $tp(N^\alpha_a /A)$ does not fork over
$A^{\beta}_1$ (by monotonicity of nonforking). Hence if $\delta$
is a limit ordinal, $tp(N^\delta_a/A)$ does not fork over
$A^\delta_1$.

But if the cofinality is smaller the union may not preserve
$\kappa$-saturation. In this case,  let   $\langle A'_1, M'_a,
M'_b,N'_a,g' \rangle$ denote the union of the respective chains;
each has cardinality $\kappa$.
  Choose $A_1 \subseteq A$ with $|A_1| = \kappa$ and such that
$A_1$ is relatively $\kappa$-saturated in $A$ and
 $A_1$ contains $A'_1$.  Then let the bound be
$\langle A_1, M_a, M_b, N_a, g\rangle$ where $M_a$ is
$\kappa$-prime over $M'_a\cup A_1$, $M_b$ is $\kappa$-prime over
$M'_b\cup A_1$, $g$ is the induced isomorphism extending $g'$
 and $N_a$ is any $\kappa$-saturated elementary extension
of $M_a \cup N'_a$ in $M$ with $N_a$ independent from $A$ over
$A_1$.
%$(N_a,A\cap N_a) \prec (M,A)$.
$\qed_{\ref{chains}}$

\begin{lemma}
\label{hammer} If $t = \langle A_1, M_a, M_b, N_a, g\rangle \in
\bK^i_{a,b}$ and $p\in S(M_a)$ is non-algebraic, orthogonal to $A$
and $p \not \perp \tp(N_a/M_a)$, then there is $t'= \langle A'_1,
M'_a, M'_b,N'_a, g' \rangle\in \bK^i_{a,b}$ with $t'$ extending
$t$ and $\tp(N_a/M'_a)$ forking over $M_a$.
\end{lemma}
Proof.  Since $M$ is $\kappa^+$-saturated,  we can find $d \in M$
realizing $p$ such that $\tp(d/N_a)$ forks over $M_a$ and $d' \in
M$ realizing $g(p)$. Now, construct $t'$ by letting $A'_1 = A_1$,
$M'_a$ be $\kappa$-prime over $M_a \cup \{d\}$, $M'_b$ be
$\kappa$-prime over $M_b \cup \{d'\}$, $g'$ be an extension of $g$
taking $d$ to $d'$, and $N'_a\prec M$ any $\kappa$-saturated
extension of $M'_a \cup N_a$. We need to show that $M'_a$ and
$M'_b$ are independent from A over $A'_1$. For this, note that
since $p\in S(M_a)$ is orthogonal to $A$ ({\em a fortiori} to
$A_1$)
 and $A$ is independent from $M_a$
over $A_1$, $d$ is independent from $A$ over $M_a$.  Since $M'_a$
is $\kappa$-prime over $M_a \cup \{d\}$, it follows that $M'_a$ is
independent from A over $A'_1$. An analogous argument shows $M'_b$
is independent from A over $A'_1$. Since $d \in M'_a$, we have
fulfilled the lemma. $\qed_{\ref{hammer}}$

For any  ordinal $\mu$ and any sequence $\langle
\abar_i:i<\mu\rangle$ and any finite $w \subseteq \mu$, $\abar_w$
denotes $\langle \abar_i:i\in w\rangle$.  We require one further
technical notion.

\begin{definition} \label{afull} We say $M_a$ is $A$-full in $M$ if for any $N$
$\kappa$-prime over $M_aA$ and for any $C_0 \subseteq M_a$,
$|C_0|\leq |T|$,
 $C_1 \subseteq A$ with $|C_1| \leq |T|$, and $C_2$ with $C_0 \subseteq C_2$,
  $C_1 \subseteq C_2 \subseteq N$, and $|C_2|\leq |T|$, there
is an elementary map $f$ taking $C_1C_2$ into $M_a$ over $C_0$
with $f(C_1) \subseteq A$ and if $C_2$ is independent from $A$
over $C_1$ then $f(C_2)$ is independent from $A$ over $f(C_1)$.
\end{definition}

We prove a characterization of a weakly benign pair; a similar
result for benign (using $\bK^1$ instead of $\bK^2$ also holds. In
view of the counterexample in given in the introduction, weakly
benign is the interesting case.
% and note that
%a variant of the proof (replacing $\bK_{a,b}$ by $\bK'_{a,b}$)
%yields the analogous result for weakly benign.

\begin{lemma}
\label{mainequiv}  Use the notation of \ref{nota}.  Suppose
 $(M,A)$ is $\kappa^+$-saturated where
$\kappa = \kappa^{|T|}$ is regular and $T=\th(M)$ is stable. The
following are equivalent.
\begin{enumerate}
\item $(M,A)$ is not weakly benign.
\item There exist $A_*, M_a, N_a, M_b,g$ contained in $M$
with $a \in M_a$, $b\in M_b$ such that:
\begin{enumerate}
\item $\langle A_*, M_a, M_b, N_a,
g\rangle \in \bK^2_{a,b}$ and $M_a \not= N_a$.
\item $M_a$ is $A$-full in $M$.
\item $N_a$ is independent from $A$ over $A_*$.
\item  $\tp(N_a/M_a)$ is orthogonal to every nonalgebraic type in $S(M_a)$ which
is orthogonal to $A$.
\item If $\dbar \in N_a-M_a$, there is no $\dbar' \in M$ which realizes
$g(\tp(\dbar/M_a))$ and such that $\dbar'$ is independent from
 $A$ over $M_b$.
 \item $M_a$ and $M_b$ are isomorphic over $A_*$ by a map $g$
 taking $a$ to $b$ and preserving strong types over $A$, i.e. $g \restriction (A^*)^{\eq}$ is the
 identity.
%\item the set $A_*$ is relatively $\kappa$-saturated inside $A$.
\end{enumerate}
\end{enumerate}
\end{lemma}
%\sidebar{ problematic if we do not assume that the  set $A$ is closed,
%best is to rep%%lace benign by weakly benign then we don't have
%a problem in this section.  from a later note from ss.
%This is addressed in Corollary 8 just below.}
%\sidebar{and to add to part 2 in the lemma
%e) the set $A_*$ is relatively $\kappa$-saturated inside $A$
%this require rereading carefully all the manuscript
%though not any non trviail change in the proof
%  ss
%Added but still need to check jtb}
%\sidebar{I am missing a big point.  Where does this argument use
%that A is closed.  I will have to modify the argument using strong types
%if I replace benign by weakly benign but I don't see where.  (This
%seems like a good way to write it; I just don't see where the difficulty
%must be fixed.)}
Note that, by general properties of orthogonality, we could
rephrase item c) as: $\tp(N_a/M_a)$ is orthogonal to every type in
$S(M_a)$ which is orthogonal to $A_*$.
%or even that $\tp(N_a/M_a)$ is orthogonal to
%every  stationary type. which is orthogonal to $A_*$.
%Before giving the proof we note one extension.
%\begin{lemma}%
%\label{mainequivcor}  Suppose $T$ is stable and $(M,A)$ is
%$\kappa^+$-saturated where $|L| \leq \kappa = \kappa^{|T|}$ is
%regular. Fix $a,b \in M-A$ and use notation of \ref{nota}.  The
%following are equivalent.
%\begin{enumerate}
%\item $(M,A)$ is not weakly benign, witnessed by $a,b$.
%\item There exist $A_*, M_a, N_a, M_b,g$ contained in $M$
%satisfying the conditions of Lemma~\ref{mainequiv} but in
 %addition $\langle A_1, M_a, M_b, N_a,
%g\rangle \in \bK'_{a,b}$.
%
%\end{enumerate}
%\end{lemma}

%Here is the proof for `benign'; for weakly benign just modify the
%proof, replacing $\bK_{a,b}$ by $\bK'_{a,b}$).

%\sidebar{looks more like proof for weakly benign???}

 Proof of
Lemma~\ref{mainequiv}: First we show that condition 2) implies
condition 1). By condition 2a),
%the definition of $\bK_{a,b}$,
 there is an $a'$
in $N_a-M_a$. Note that since $A^*$ is relatively
$\kappa^+$-saturated in $A$ and $M_a$ ($M_b$) is independent from
$A$ over $A^*$, $M_a \cap A  =M_b \cap A=A^*$.    It follows that
$(g \cup \id)\restriction \acl(A^{eq})$ is an elementary map in
$L^{\eq}$. Let $\abar=\langle a_i:i<\kappa\rangle$ enumerate
$M_a-A$ with $a_0 = a$; denote $g(a_i)$ by $b_i$ so $\bbar=\langle
b_i:i<\kappa\rangle$ enumerates $M_b$.   For any finite set of
$L$-formulas $\Delta$ and finite subset $w$ of $\kappa$, let
$\phi_{\Delta,w}(\xbar;a',\abar_w,\bbar_w)$ be the $L^*$-formula
which assert that $x \bbar_w$ and $a'\abar_w$ realize the same
$\Delta$-type over $A$.  %Since $g$ preserves strong types and
%$M_a$, $M_b$ is each independent from $A$ over $A^*$,
For any finite $w$, $\abar_w$ and $\bbar_w$ realize the same
$L$-type over $A$.

 Now, let $q
=\{\phi_{\Delta,w}(\xbar;a',\abar_w,\bbar_w):0\in w \subset_\omega
\kappa, \Delta \subset_\omega L\}$. Putting $0\in w$ guarantees
$a,b$ are in any relevant $\abar_w,\bbar_w$.
  So $q $ is a set
of $\kappa$ $L^*$-formulas with free variable $x$ and parameters
from $M_a \cup M_b\cup \{ a'\}$. If $q$ is finitely satisfied in
$(M,A)$, then $q$ is realized in $M$ by some $b'$, since $M$ is
$\kappa$-saturated as an $L^*$-structure. But since $a'$ is
independent from $A$ over $M_a$,  $b'$ realizes the unique
nonforking extension of $g(\tp(a'/M_a))$ to $M_b\cup A$
contradicting condition d). If $q$ is not finitely satisfiable,
there is a formula $\phi_{\Delta,w}$ which demonstrates the $L^*$
type of $\abar_w$ and $\bbar_w$  over $A$ are different.
%This implies $a$ and $b$ realize different $*$-types over $A$???

%there is an $L$-formula $\psi\in \Delta$ which witnesses that
%there is no $b'\in M$ with $\abar a' \equiv_A \bbar b'$, i.e.
%$\models \psi(\abar a')$ but there is no $b'$ with $\models
%\psi(\bbar b')$.

To show the converse, we suppose that $\abar$ and $\bbar$ realize
the same (strong)-type over $A$ but that there is an $a'$ such
that there is no $b'\in M$ with $\abar a' \equiv_{A,L} \bbar b'$.

We will use the following basic fact:
\begin{fact}
\label{fact}
\begin{enumerate}
\item If $A_1$ is relatively $\kappa$-saturated in $A$ and
$C$ is independent from $A$ over $A_1$, then
$CA_1$ is relatively $\kappa$-saturated in $CA$.
\item If $A_1$ is relatively $\kappa$-saturated in $A$ and
$D$ is $\kappa$-atomic over $A_1$,
$D$ is independent from $A$ over $A_1$.
\end{enumerate}
\end{fact}

%\sidebar{Whoops, I overestimated my understanding.  Following is a
%pretty literal description of a page labeled F529 in center of
%page and on right 2002/1013.

%It is supposed to say how to alter the argument to guarantee the
%property I have labeled A-full.  I suppose that it should be put
%in as part of the lemma below but I don't see it now.  I also
%think there will have to be a limit stage argument paragraph i.e
%in before 1.4.  It would be most elegant to add $A$-full as part
%of the definition of $K_{a,b}$.

%Here is your proof.

%Given $t_\alpha$, choose $N_\alpha \prec M$, $\kappa$-prime over
%$M^\alpha_aA$ (unique up to isomorphism). Choose $A_*^{\alpha+1}
%M_a^{\alpha+1} M_b^{\alpha+1}$ such that $A_*^{\alpha} \subseteq
%A_*^{\alpha+1}\subseteq A$, $|A_*^{\alpha+1}| = \kappa$ and closed
%enough, $ M_a^{\alpha+1}$ is $\kappa$-prime over $
%M_a^{\alpha}A_*^{\alpha+1}$, $ M_b^{\alpha+1}$ is $\kappa$-prime
%over $ M_b^{\alpha}A_*^{\alpha+1}$

%and by bookkeeping satisfies a case of (*)  of or just $
%M_a^{\alpha}A_*^{\alpha} \subseteq  (M_a^{\alpha},A_*^{\alpha+1})
%\prec_{L_(|T|^+,|T|^+} (N^\alpha,A)$ (as $\kappa = \kappa^|T|$
%possible and can add being prime $ M_a^{\alpha+1}$ is
%$\kappa$-prime over $ M_a^{\alpha}A_*^{\alpha+1}$.
% \underline{The point} No problem to choose $
%M_b^{\alpha+1}$  such that $g_{\alpha +1}$ exists because of
%$\kappa$-primeness. }

The following lemma essentially shows 1) implies 2) of
Lemma~\ref{mainequiv}.

\begin{lemma}\label{crux}  There is  a $t = \langle A_*, M_a, M_b, N_a,
g\rangle \in \bK^2_{a,b}$ such that
\begin{itemize}
%\item $N_a$ is independent from $A$ over $A_*$.
\item[A]
$N_a \not= M_a$,
\item [B] $\tp(N_a/M_a)$ is orthogonal to every nonalgebraic type
 in $S(M_a)$,
which is  orthogonal to
$A$.
\item  [C] If $\dbar \in N_a-M_a$, there is no $\dbar' \in M$ which realizes
$g(\tp(\dbar/M_a))$ and such that $\dbar'$ is independent from
 $A$ over $M_b$.
\item[D]
$M_a$ is $A$-full.
\end{itemize}
\end{lemma}
Proof.
 Try to construct by induction a sequence
 $\langle t_\alpha:\alpha
< \kappa^+\rangle$ where  $t_\alpha = \langle A^\alpha_*,
M^\alpha_a, M^\alpha_b, N^\alpha_a, g^\alpha\rangle $ of elements
of $\bK^2_{a,b}$ which are increasing in the natural partial
order, continuous at limit ordinals of cofinality greater than
$\kappa_r(T)$ and with $a' \in N^0_a$.

\begin{enumerate}
\item
If $\alpha$ is an even ordinal there are several cases.

\begin{enumerate}
\item Suppose condition B fails, i.e.
 for some $d\in N_a$, $p = \tp(d/M_a)$
is nonorthogonal to some stationary type $q\in S(M_a)$ which is
orthogonal to $A$.
 Then by Lemma~\ref{hammer},
there is       $t'= \langle A'_*, M'_a, M'_b,N'_a, g' \rangle\in
\bK^2_{a,b}$ with $t'$ extending $t$ and $\tp(N_a/M'_a)$ forks
over $M_a$.
\item Suppose condition B holds.
\begin{enumerate}
\item If
$\alpha$ is a limit ordinal of cofinality $\kappa$, stop. \item If
$\alpha$ is a limit ordinal of cofinality $<\kappa$  or $\alpha$
is a successor ordinal, let $t_{\alpha+1} = t_{\alpha}$.
\end{enumerate}
\end{enumerate}
\item $\alpha$ is an odd successor ordinal.  Choose  an auxiliary
${\hat M_a}^{\alpha}$
%\sidebar{no cardinality requirement on $N^\alpha$; $N^\alpha$
%already chosen!! }
$\kappa$-prime over $M^\alpha_a A$.  Choose $A_*^{\alpha+1},
M_a^{\alpha+1}, M_b^{\alpha+1}$ such that $A_*^{\alpha} \subseteq
A_*^{\alpha+1}\subseteq A$, $|A_*^{\alpha+1}| = \kappa$ and so
that $$({ M_a}^{\alpha+1},A_*^{\alpha+1}) \prec_{L_(|T|^+,|T|^+}
({\hat M_a}^{\alpha},A)$$ and
 $ M_a^{\alpha+1}$ is $\kappa$-prime over $
M_a^{\alpha}A_*^{\alpha+1}$.  This is possible since $\kappa =
\kappa^{|T|}$.  In particular, $ M_a^{\alpha+1}$ is independent
from $A$ over $A_*^{\alpha+1}$. The $\kappa$-primeness allows us
to easily construct $M_b^{\alpha+1}$ and $g_{\alpha+1}$.  Now
choose $N_a^{\alpha+1}$ to be a $\kappa$-saturated extension of $
M_a^{\alpha+1}$ that is independent from $A$ over $A_*^{\alpha+1}$
\item If $\alpha$ is a limit ordinal choose $t_\alpha$ by
Lemma~\ref{chains}.
\end{enumerate}

% and satisfies the following three
%conditions:
% \begin{itemize}
% \item [A.]If $\cf(\alpha) = \kappa$, $M^{\alpha+1}_a$ depends on
% $N^{\alpha}_a$ over $M^{\alpha}_a$;
% \item [B.]
% otherwise: $\tp(N^{\alpha}_a/A)$ does not fork over
% $A^{\alpha+1}_*$.

% \item [C.] $a' \in N^0_a$.
% \end{itemize}

We cannot carry out this construction for $\kappa^+$ steps. If we
did,  by clause 1) of the construction at each limit $\alpha$ with
$\cf(\alpha) = \kappa$,  clause B) fails.  Thus, $M^{\alpha+1}_a$
depends on
 $N^{\alpha}_a$ over $M^{\alpha}_a$ for all such
 $\alpha$, which contradicts stability.
  (If we were dealing with finite sequences, the
bound would be $\kappa(T)$; since we deal with sets of cardinality
$\kappa$, the bound is $\kappa^+$.)

 Fix $\alpha$ where the construction stops. We have constructed $t_\alpha =
\langle A^\alpha_*, M^\alpha_a, M^\alpha_b, N^\alpha_a,
g^\alpha\rangle$ but for any choice of $t_{\alpha+1} \in
\bK^2_{a,b}$, $M^{\alpha+1}_a$ is independent from $N^{\alpha}_a$
over $M^{\alpha}_a$.  Note that each member of $t_\alpha = \langle
A^\alpha_*, M^\alpha_a, M^\alpha_b, N^\alpha_a, g^\alpha\rangle$
is the union of the respective member of $t_\beta$ over $\beta<
\alpha$. We claim this $t_\alpha$ is a $t$ satisfying the
conditions of the lemma.

 %when we show $t_\alpha \in \bK'_{a,b}$.
%We check that  $ N^\alpha_a$, $  A$ are
%independent
% over  $ A^\alpha_*$:  By clause (B) for every
%$\alpha<\beta$, $tp(N^\alpha_a /A)$ does not fork over
%$A^{\beta}_*$ (by monotonicity of nonforking). Hence if $\delta$
%is a limit ordinal, $tp(N^\delta_a/A)$ does not fork over
%$%A^\delta_1$.  In particular, $ N^\alpha_a$, $  A$ are independent
 %over  $ A^\alpha_1$.

For clause A note
 $N^{\alpha}_a \neq M^{\alpha}_a$ since
$a' \in N^{\alpha}_a$ and $a'$ cannot be in the domain of
$g^\alpha$ by the original choice of $a'$.  Since the construction
stopped clause B, holds.

For clause C,  we must show that if $\dbar \in N_a-M_a$, there is
no $\dbar' \in M$ which realizes $g(\tp(\dbar/M_a))$ and such that
$\dbar'$ is independent from
 $A$ over $M_b$. Fix $\dbar \in N_a-M_a$; if such a $\dbar'$ exists, choose $M^{\alpha+1}_a,
M^{\alpha+1}_b$ contained in $M$ prime over $M^{\alpha}_a \dbar$ and
 $M^{\alpha}_b \dbar'$  respectively.  We easily extend $g^\alpha$ to
$g^{\alpha+1}$ mapping $M^{\alpha+1}_a$ to $ M^{\alpha+1}_b$. By
the construction,
 $ A^\alpha_*$ is relatively $\kappa$-saturated in $ A$.
So, $ M^\alpha_a\cup \{\dbar\}$ and $ A$ are independent over $
A^\alpha_*$ by monotonicity, as $N^\alpha_a$ is independent from
$A$ over $A^\alpha_*$. Now by Fact~\ref{fact} 1),
 $ M^\alpha_a \cup \{\dbar\}$  is relatively $ \kappa$ -saturated
inside $ M^\alpha_a \cup \{\dbar\}  \cup A$. Whence,  by
Fact~\ref{fact} 2) $ M^{\alpha+1}_a$  and  $ A$ are independent
over $ M^\alpha_a\cup \{\dbar\} $.  By transitivity of nonforking,
$ M^{\alpha+1}_a$ and $ A$  are independent over $ A^{\alpha}_*$.
Similarly, since  $\dbar'$ is independent from
 $A$ over $M_b$, $M^{\alpha+1}_b$ is
independent from $A$ over $A^{\alpha} _*$.  But now, $\dbar \in
(M^{\alpha+1}_a \cap N^\alpha_a) - M^\alpha_a$ $so N^\alpha_a$
depends on $M^{\alpha+1}_a$ over $M^{\alpha}_a$ and we have
violated the choice of $\alpha$.

Finally we verify clause D: $M_a$ is $A$-full.  Choose $N$, which
is $\kappa$-prime over $AM_a$.  Then $N$ can be embedded over
$AM_a$ into ${\hat M}^\alpha_a=\cup_{i< \alpha}{\hat M}^i_a$. By
the Tarski union of chains theorem (using clause 2) of the
construction), $(M^\alpha_a, A \cap M^\alpha_a)
\prec_{L_{|T|^+,|T|^+}} ({\hat M}^\alpha_a,A)$. Let $C_0, C_1, C_2
\subseteq N$ satisfy the hypotheses of the definition of $A$-full.
The elementary submodel condition easily allows us to define the
required function $f$. $\qed_{\ref{crux}}$

\section{The Superstable Case}

The aim of  this section is to prove that if $M$ is a model of a
superstable theory and $A \subset M$, then $(M,A)$ is weakly
benign.  This is a
%modest
generalization of a result of Bouscaren \cite{Bouscarenpairs}, who
showed, in our terminology that every {\em submodel} of a
superstable structure is benign.
%established
%the conclusion with a different proof under the additional
%hypothesis that $A$ is the universe of a model.

\begin{theorem}  If $M$ is a model of
a superstable theory and $A \subset M$, then $(M,A)$ is weakly
benign.
\end{theorem}

Proof.   We work in $\Mscr^{\eq}$.  Without loss of generality,
assume $(M,A)$ is $\kappa^+$-saturated for a regular $\kappa$
satisfying $\kappa^{|T|} =\kappa$.  By Lemma~\ref{mainequiv} if
$(M,A)$ is not weakly benign, there exist $A_*, M_a, N_a, M_b,g$
contained in $M$ satisfying the conditions of
Lemma~\ref{mainequiv} and with $\langle A_*, M_a, M_b, N_a,
g\rangle \in \bK^2_{a,b}$.

Since $M_a$ is properly contained in $N_a$, we can choose $\cbar
\in M_a$ and  $\phi(x,\cbar)$ to have minimal $D$-rank among all
formulas with $\phi(N_a,\cbar) \neq \phi(M_a,\cbar)$.  Then for
any $d^* \in \phi(N_a,\cbar)\setminus \phi(M_a,\cbar)$,
$p^*=\tp(d^*/M_a)$ is  regular. Without loss of generality again,
we can fix $d^*$, which does not fork over $\cbar$ and so that
$p^*$ has the same $D$-rank as $\phi(x,\cbar)$ and
$\tp(d^*/\cbar)$ is stationary. By clause c) of
Lemma~\ref{mainequiv}, $p^*$ is not orthogonal to $A_*$.  So,
there is a $q' \in S(M_a)$ which  does not fork over $A_*$ and is
nonorthogonal and so non-weakly orthogonal to $p^*$. Fix $C
\subseteq A_*$ with $|C| \leq |T|$ and $\cbar$ is independent from
$A_*$ over $C$.
 Without loss of generality $\tp(\dbar^*/A_*\cbar) \not \perp^w q \restriction
(A_*\cbar)$ and $\tp(\dbar^*/C\cbar) \not \perp^w q\restriction
(C\cbar)$.   Let $\Pscr = \{p:  p  {\rm \ is\ regular,\
stationary,\ and\ nonorthogonal\ to \ }  p^*   \}$. $\Pscr$ is
based on $B = \acl^{\eq}(C)$, i.e. every automorphism of $\Mscr$
fixing $B$ maps $\Pscr$ to itself.

%\sidebar{Alternatively to conform with book, let $\Pscr = \{p: p
%is regular, stationary, and nonorthogonal to p^* \}$. $\Pscr =
%\{F(p^*): F \in \aut_{C^{\eq}}(\Mscr^{\eq})\}$. }

 If $\cbar''\in M$
realizes $\tp(\cbar/A^{\eq})$ and $d''\cbar'' $ realizes $r =
\tp(d^*\cbar/B)$, then $\tp(d''/\cbar'')$ is regular and
nonorthogonal to $p^*$.  We can find $\langle \cbar_i:i<
\omega\rangle$ in $M_a$ with $\cbar_0 =\cbar$ which are
indiscernible over $B$ and which are based on $B$. The
$r(\xbar,\cbar_i)$ are regular, pairwise nonorthogonal, and all
nonorthogonal to $\Pscr$ and each $r(\xbar,\cbar_i)$ is not weakly
orthogonal to $q'\restriction (B\cbar_i)$. Note
$r(\xbar,\cbar_i)\subset p^*$.  Let $r_i \in S(M)$ denote the
nonforking extension of $r(\xbar,\cbar_i)$ to
 $S(M)$.
% Without
%loss of generality we may assume the nonforking extensions of
%$r(\xbar,\cbar_i), r(\xbar,\cbar_j)$ to $ S(C\cbar_i\cbar_j)$ are
%not weakly orthogonal. (Choose $\dbar_{i,j}$ for $i,j < \omega$
%with $\dbar_{i,j}$ realizing $r(\xbar,\cbar_i)$ such that for each
%$j$,
%$$\tp(\dbar_{i,j}/\bigcup_{\ell<\omega}\cbar_\ell \cup \{\dbar_{i',j'}:(i'<
%\ell \wedge j'< \omega)\vee (i'=i \wedge j'<j\} \cup C)$$ does not
%fork over $\cbar_j$.  For large enough $m$ replace $\cbar_i$ by
%$\cbar_i\dbar_{i,0} \ldots \dbar_{i,m})$.
 By Section V.4 of
\cite{Shelahbook2nd}, there is a $q \in S(B)$, which
 is $\Pscr$-simple and $k < \omega$ such that $w_{\Pscr}(q) > 0$ and $q(\Mscr) \subseteq
 \acl(B \cup \bigcup_{i<k}\cbar_i \cup
 \bigcup_{i<k}r(\Mscr,\cbar_i)$. (This $q$ is actually $q'/E$ for an appropriate
definable (over $B$) equivalence relation; compare V.4.17(8) of
\cite{Shelahbook2nd}.)

 Let $q^+$ denote the unique nonforking extension of $q$ to
 $S(M)$,
$p_a^+$ denote the unique nonforking extension of $p^*$ to
 $S(M)$, and
$p_b^+$ denote the unique nonforking extension of $g(p^*)$ to
 $S(M)$.
Clearly, $p_a^+\restriction (M_a \cup A)$ is a nonforking
extension of the stationary type $p^*$ and realized by $\dbar^*$;
so it is equivalent to $p_a^+\restriction \acl(M_a \cup A)$.

\begin{remark}\label{anvil}Note $(g \cup \id_A)(p_a^+\restriction (M_a \cup A)=
p_b^+\restriction (M_b \cup A) \sim p_b^+\restriction \acl(M_b
\cup A)$ is omitted in $M$.
\end{remark}

 We use the next lemma several times.

\begin{lemma}\label{d4Balpha}
If $A^{\eq}\subseteq N_1 \subseteq N_2 \subseteq M$ and $N_1, N_2$
are $|T|^+$-saturated then
$$w_{\Pscr}(q(N_2),N_1) =w_{\Pscr}(q(N_2),q(N_1)A^{\eq}).$$
\end{lemma}

Proof. Fix $\bbar \in N_1$ and choose $D \subseteq q(N_1)A^{\eq}$
with
 $|D| \leq |T|$ such that $\tp(\bbar/q(N_1)A^{\eq})$ does not fork
 over $D$.  If $\tp(\bbar/q(N_2)A^{\eq})$
 forks over $D$, there are finite $\dbar_1 \subseteq q(N_2)$ and $\dbar_2
 \subseteq A^{\eq}$ such that $\tp(\bbar/B D \dbar_1\dbar_2)$ forks
 over $D$.  But there is a $\dbar' \in q(N_1)$ realizing
 $\stp(\dbar_1/D\bbar\dbar_2)$, which contradicts $\tp(\bbar/q(N_1)A^{\eq})$ does not fork
 over $D$.

So $\tp(\bbar/q(N_2)A^{\eq})$ does not fork over $q(N_1)A^{\eq}$.
Since $\bbar$ was arbitrary in $N_1$, $\tp(N_1/q(N_2)A^{\eq})$
does not over $q(N_1)A^{\eq}$. By symmetry of forking,
$\tp(q(N_2)/N_1A^{\eq})$ does not fork
 over $q(N_1)A^{\eq}$.  Since $A^{\eq} \subseteq N_1$ we finish.
 $\qed_{\ref{d4Balpha}}$

The proof now proceeds by a series of claims.  The key idea is
that $w_{\Pscr}(q(M), A^{\eq})$ can be calculated as either
$w_{\Pscr}(q(M),q(M_b)\cup A^{\eq})+ w_{\Pscr}(q(M_b), A^{\eq})$
or as
 $w_{\Pscr}(q(M),q(M_a)\cup A^{\eq})+ w_{\Pscr}(q(M_a), A^{\eq})$. We will
calculate both ways to obtain a contradiction.  We begin with the
$M_a$ side.

\begin{claim}\label{Hebrew}
If
 $\dim (r_0\restriction A_* \cbar_0 , M_a)$ is finite, then
$w_{\Pscr}( q(M_a), A_* \cup \bigcup_{i<k}\cbar_i)$ is finite.
\end{claim}

Proof.  If $u$ is a finite subset of $\omega$, since the $r_i$ are
regular, it is easy to show
 that for each $i$,  $\dim(r_i \restriction
 (A_* \cbar_i),M_a)$ is finite iff $\dim(r_i\restriction
 (A_* \cup \cbar_i \cup_{ j \in u}\cbar_j),M_a)$ is finite.
Since the $r_i \restriction (A_* \cbar_i \cbar_j)$ are regular
 and pairwise not weakly orthogonal
 $$\dim (r_i\restriction A_* \cbar_i \cbar_j, M_a)
 =\dim (r_j\restriction A_* \cbar_i \cbar_j, M_a).$$
The previous two sentences imply:
 $\dim (r_i\restriction A_* \cbar_i , M_a)$ is finite iff
  $\dim (r_j\restriction A_* \cbar_j , M_a)$ is finite.
 So if
 $\dim (r_0\restriction A_* \cbar_0 , M_a)$ is finite then
 $w_{\Pscr}(
 \bigcup_{i<k}r_i(M_a,\cbar_i), A_* \cup \bigcup_{i<k}\cbar_i)$ is
 finite; whence $w_{\Pscr}(
q(M_a), A_* \cup \bigcup_{i<k}\cbar_i)$ is finite.
$\qed_{\ref{Hebrew}}$

Now we drop the $\bigcup_{i<k}\cbar_i$ in the conclusion.

%{following from aleph through dalet}
\begin{claim}\label{i}
 $\dim (r_0\restriction A_* \cbar_0 , M_a)$ is finite implies
$w_{\Pscr}(q(M_a), A_*)$ is finite. \end{claim}
 Proof. Find $\dbar \subseteq q(M)$ such that $
\bigcup_{i<k}\cbar_i$ is independent from $A_* \cup q(M)$ over
$A_*\cup \dbar$.  Now, as $\tp(\dbar/A_*)$ is $\Pscr$-simple,
$w_{\Pscr}( q(M_a), A_*)=w_{\Pscr}(q(M_a),
A_*\dbar)+w_{\Pscr}(\dbar, A_*)$.  The second term is finite and
$w_{\Pscr}( q(M_a), A_*\dbar)=w_{\Pscr}( q(M_a), A_*\dbar\cup
\bigcup_{i<k}\cbar_i)$ by the independence. But, $w_{\Pscr}(
q(M_a), A_*\dbar\cup  \bigcup_{i<k}\cbar_i)= w_{\Pscr}( q(M_a),
A_*\cup \bigcup_{i<k}\cbar_i)-w_{\Pscr}(\dbar, A_*\cup
\bigcup_{i<k}\cbar_i)$.  Now the first of the last two terms is
finite by Claim~\ref{Hebrew} (since $\dim (r_0\restriction A_*
\cbar_0 , M_a)$ is finite) and the second by the finiteness of
$\dbar$ so $w_{\Pscr}(q(M_a), A_*)$ is finite . $\qed_{\ref{i}}$

\begin{claim}\label{h} $\dim(r_0,M_a)$ is finite.
\end{claim}
Note that $p^+_a\restriction(B\cbar_0) =
r_0\restriction(B\cbar_0)$.
 Choose by induction
$\abar_\alpha \in M_a $ so that $\abar_\alpha$ realizes
$p^+_a\restriction A^{\eq}_* \cup g(\cbar_0) \cup
\{\abar_\beta:\beta< \alpha\}$ for as long as possible to
construct:  ${\bf I} = \langle \abar_\alpha:\alpha <
\alpha^*\rangle$. Clearly $\alpha^* <|M_a|^+$, but in fact
$\alpha^*$ is finite. As, since $M_a$ is independent from $A$ over
$A_*$, ${\bf I}$ is a set of indiscernibles over $A$. Since $M$ is
$\kappa^+$-saturated, if ${\bf I}$ is infinite $\langle
g(\abar_\alpha):\alpha < \alpha^*\rangle$ can be extended to a set
${\bf J}$ of indiscernibles over $A$ contained in $M_b$ with
cardinality $\kappa^+$. Then all but at most $\kappa$ members of
 ${\bf J}$  realize $p_b^+\restriction (M_b \cup A)$ contradicting
 Remark~\ref{anvil} that $p_b^+\restriction (M_b \cup A)$ is omitted
 in $M$.
 $\qed_{\ref{h}}$

Now, easily we have
\begin{claim}\label{claimnewj}
$w_{\Pscr}(q(M_a), A_*)=w_{\Pscr}(q(M_a), A^{\eq})$ is finite.
\end{claim}
The equality holds by the independence of $M_a$ and $A$ over
$A_*$. The finiteness follows from Claim~\ref{h} and
Claim~\ref{i}.
 $\qed_{\ref{claimnewj}}$

%clauss ell

%$$w_{\Pscr}(q(M_b), A^{\eq})= w_{\Pscr}(q(M_b), A_*) < \omega.$$
%The equality holds by the independence of $M_b$ and $A$ over $A_*$
%and the inequality by Claim  \ref{claimnewj}.

%\sidebar{clause m}

The next claim involves both $M_a$ and $M_b$.

\begin{claim}\label{hey}  Suppose $w_{\Pscr}(q(M_a),A_*)$ is
finite and $N \prec M$ is $\kappa$-prime over $M_b  A$. \\ Then
$w_{\Pscr}(q(N), q(M_b) A)=0$.
\end{claim}

%\sidebar{literal copying - can't interpret

%is $B^*$ earlier B?}
Proof. Since $w_{\Pscr}(q(M_a),A_*)$ is finite, and $A$, $M_a$ are
independent over $A_*$,we can choose finite $D\subseteq q(M_a)$
with $w_{\Pscr}(q(M_a),A_*)
=w_{\Pscr}(q(M_a),A)=w_{\Pscr}(D,A_*)=w_{\Pscr}(D,A)$.

Now assume for contradiction that $w_{\Pscr}(q(N), q(M_b) A)>0$.
Let $N' \prec M$ be $\kappa$-prime over $M_a \cup A$, so there is
$g^+ \supseteq g \cup {\rm id}_A$ which is an isomorphism from
$N'$ onto $N$.  Then there is a finite $D_2 \subseteq q(N')$ with
$w_{\Pscr}(D_2, M_a A)>0$.
%and since $A$, $M_a$ are independent over $M_aA_*$,
Choose $C_0 \subseteq M_a$, $|C_0|\leq |T|$ with $DB \subseteq C_0$ and
 $C_1 \subseteq A$ with $|C_1| \leq |T|$ so that $D_2$ is independent
from $M_aA$ over $C_0C_1$ and is the unique nonforking extension
of $\tp(D_2/C_0C_1)$ to $S(M_aA)$ which is realized in $M$. Recall
that $M_a$ is $A$-full and apply the Definition~\ref{afull} of
$A$-full with $C_0C_1D_2$ playing the role of  $C_2$ to obtain an
embedding $f$. Then, $f(D_2) \subseteq q(M_a)$ and $f(D_2)$ is
independent from $C_0A$ over $C_0f(C_1)$. Thus,
$$w_{\Pscr}(f(D_2),AD)=w_{\Pscr}(D_2,AD)\geq w_{\Pscr}(D_2,q(M_a)A)>0.$$
This implies $w_{\Pscr}(q(M_a),A)\geq w_{\Pscr}(Df(D_2),A)=
w_{\Pscr}(D,A)+w_{\Pscr}(f(D_2),AD)> w_{\Pscr}(D,A)$, which
contradicts our original choice of $D$.
 $\qed_{\ref{hey}}$

%So it suffices to So we have condition $AA$.  $q
%\in S(B^*)$ is stationary, $B^* \subseteq A_* \subseteq M_a$, $q$
%is $\Pscr$-simple, $\Pscr$ based on $B^*$,  $w_{\Pscr}(q(M), A_*)$
%is finite and there is an $N' \prec N$ be $\kappa$-prime over $M_a
%\cup A$ with $w_{\Pscr}(q(N'), q (M_b) A)>0$
%\sidebar{last q makes no sense or should there be parentheses?}
% and we choose a finite
%$D \subseteq q(M_a)$ such that $w_{\Pscr}(q(M_a),
%A_*)=w_{\Pscr}(D, A_*)$. This contradicts the new clause added to
%the definition in section 1. $\qed_{\ref{hey}}$

%Now we do the calculation on the $M_a$ side. probably wrong

\begin{claim}\label{e}
$w_{\Pscr}(q(M), q(M_b) A)=0$
\end{claim}
Let $N \prec M$ be $\kappa$-prime over $M_b \cup A$, so
$p_b^+\restriction (M_b \cup A)$ has a unique extension in $S(N)$.
If $w_{\Pscr}(q(M),N)>0$ then for some $\bbar \in q(M)$,
$w_{\Pscr}(\bbar,N)>0$ so $\tp(\bbar/N) \not \perp p^+_b$; recall
$p^+_b$ is parallel to $p^+_b\restriction N$.  So
$p^+_b\restriction N$ is realized in $M_b$ contradicting
Remark~\ref{anvil}. Now $0 =w_{\Pscr}(q(M),N)$ which equals
$w_{\Pscr}(q(M),q(N)A^{\eq})$ by Lemma~\ref{d4Balpha}. Since
$A^{\eq} \subseteq N_b \subseteq N \subseteq M$,
$$w_{\Pscr}(q(M),q(M_b)A^{\eq})=w_{\Pscr}(q(M),q(N)A^{\eq})+w_{\Pscr}(q(N),q(M_b)A^{\eq})=
0+0=0.$$  The first 0 was noted in the previous sentence and
the second is Claim~\ref{hey}. $\qed_{\ref{e}}$

Now calculating with respect to $M_b$, we have:

\begin{claim}\label{bside} $w_{\Pscr}(q(M), A^{\eq})= w_{\Pscr}(q(M_b), A^{\eq})$
is finite.
\end{claim}
Proof.
$$w_{\Pscr}(q(M), A^{\eq})= w_{\Pscr}(q(M),q(M_b)A^{\eq})+ w_{\Pscr}(q(M_b), A^{\eq})$$
$$=0 + w_{\Pscr}(q(M_b), A^{\eq})< \omega.$$
The first equality holds by additivity \cite{Shelahbook2nd} and
Lemma~\ref{d4Balpha}, the second by Claim~\ref{e}, and the third
by the last observation. $\qed_{\ref{bside}}$

Now we analyze using  $M_a$.
\begin{claim}\label{g} $w_{\Pscr}((q(M),q(M_a) \cup A)\geq 1$.
\end{claim}

Proof.
% \sidebar{clause f}
$w_{\Pscr}(\dbar^*, M_a \cup A)\geq 1$ since $\dbar^*$ is
independent from $A$ over $M_a$.
 Let
$N$ be $\kappa$-prime over $M_aA^{\eq}$.  As
$\tp(\dbar^*/M_aA^{\eq})$ has all its restrictions to set of size
less than $\kappa$ realized in $M_aA^{\eq}$, $\tp(\dbar^*/N)$ does
not fork over  $M_aA^{\eq}$.  Thus, $\dbar^*$ realizes
$p^+_a\restriction N$.  Since $p^+_a\restriction N$ is not
orthogonal to $q^+\restriction N$, there is $\bbar \in q^+(M)$
which depends on $\bbar$ over $N$. So $w_{\Pscr}(\bbar,N) >0$
whence $w_{\Pscr}(q(M), N)>0$. By monotonicity,
$w_{\Pscr}((q(M),q(M_a) \cup A_*)\geq w_{\Pscr}(q(M), q(N)
A^{\eq})$.  But, by Lemma~\ref{d4Balpha}, $w_{\Pscr}(q(M), q(N)
A^{\eq}) = w_{\Pscr}(q(M), N) >0$. $\qed_{\ref{g}}$

Now we have

\begin{equation} w_{\Pscr}(q(M), A^{\eq})=w_{\Pscr}(q(M),q(M_a) A^{\eq})+
w_{\Pscr}(q(M_a), A^{\eq})\geq  1  + w_{\Pscr}(q(M_a), A^{\eq})<
\omega.\label{bigeq}\end{equation} Here, the  first equality is by
\cite{Shelahbook2nd} and Lemma~\ref{d4Balpha} and the second by
Claim~\ref{g}.  The finiteness comes from Claim~\ref{claimnewj}.
Since $g \cup \id_{A^{\eq}}$ is an elementary map,
$w_{\Pscr}(q(M_a), A^{\eq})= w_{\Pscr}(q(M_b), A_*)$. We
substitute in Equation~\ref{bigeq}, using Claim~\ref{bside}:

$$w_{\Pscr}(q(M_a), A^{\eq})=w_{\Pscr}(q(M), A^{\eq})= w_{\Pscr}(q(M_a),
A^{\eq})+1,$$ or subtracting, $0=1$ so we finish.

\end{document}